\input amstex

\documentstyle{amsppt}
\loadbold
\loadeufm
\magnification=1200
\NoBlackBoxes

\topmatter
\title A Negative Answer to A Question by Rieffel
  \endtitle
\author Jae-Kwan Shim \endauthor
\address
\email jshim\@semi.korea.ac.kr \endemail
\endaddress
\abstract In this article, we address one of the questions raised
by Marc Rieffel in his collection of questions on deformation
quantization. The question is whether the
$K$-groups remain the same under flabby strict deformation
quantizations.
By ``deforming" the question
slightly, we produce a negative answer to the question.
\endabstract
\affil Department of Industrial Applied Math  \\
 Kyungpook National University \\
 Taegu, Korea
 \endaffil
 \thanks This work was supported by the Brain Korea 21 Project. \endthanks
\endtopmatter\document


\newcount\refno
\refno=0

\newbox\References
\setbox\References=\vbox{}

\def\refs#1=#2\par{\advance\refno by 1 %
   \edef#1{{\bf[\the\refno]}}
   \setbox\References=\vbox{\unvbox\References\noindent\hangindent=3pc
   \hangafter=1%
   \hbox to 20pt{\hss\bf[\the\refno]\kern 5pt}#2\unskip\nobreak\strut\par}}

\def\biblio{\nobreak\medskip\unvbox\References\par\smallskip}


\refs\Asp={P.~Aspinwall: {\sl $K3$ surfaces and string duality},
hep-th/9611137.}

\refs\CDS={A.~Connes, M.~Douglas, and A.~Schwarz: {\sl
Noncommutative geometry and matrix theory: compactification on
tori}, J.~High Energy Phys.~{\bf no 2} (1998) Paper3.}

\refs\FaW={C.~Farsi and N.~Watling: {\sl Symmetrized
non-commutative tori}, Math.~Annalen {\bf 296} (1993), 739--741.}

\refs\Kont={M.~Kontsevich: {\sl Deformation quantization of
Poisson manifolds I}, \newline q-alg/9709040.}

\refs\RaWi={I.~Raeburn and D.~Williams: {\it Morita Equivalence
and Continuous Trace $C^*$-algebras}, Mathematical Surveys and
Monographs, {\bf 60}, American Mathematical Society, Providence,
 1998.}

\refs\RiefNon={Marc A.~Rieffel: {\sl Non-commutative tori -- a case
study of non-commutative differential manifolds}, Contemp.~Math
{\bf 105} (1990), 191--211.}

\refs\RiefPro={Marc A.~Rieffel: {\sl Projective modules over
higher-dimensional non-commutative tori}, Can.~J.~Math, {\bf XL},
No.~2 (1988), 257--338.}

\refs\RiefQue={Marc A.~Rieffel: {\sl Questions on quantization,}
Contemp.~Math {\bf 228} (1998), 315--326.}

 In his collection    of questions on deformation quantization
 \RiefQue,
 Marc Rieffel asked the following;  ``Are the $K$-groups of the
 $C^*$-algebra completions of the algebras of any flabby strict
 deformation quantization all isomorphic?"  Up to my knowledge,
 the question is still open. But this article will show that the
 answer is negative if we ask the same question for the
 case of orbifolds.

\definition{Definition 1 \RiefQue}
 Let $(M, \{\cdot, \cdot \} )$ be a Poisson manifold.
 A strict deformation quantization of $M$ in the direction of
 $\{\cdot, \cdot \}$
is a dense $*$-algebra $A$ of $C^{\infty}(M)$ which is closed
under the Poisson bracket, together with a closed subset $I$ of
the real line containing 0 as a non-isolated point, and for each $\hbar\in I$
an associative product $*_\hbar$, an involution ${\ }^{*_\hbar}$,
and a pre-$C^*$-norm $\| \cdot \|_\hbar$ on $A$, which for $\hbar=0$
are the original pointwise multiplication, complex conjugation,
and supremum norm respectively, and such that \newline
(1) for each $f\in A$, $f\to \|f\|_\hbar$ and $\|f\|^{*_\hbar}$
are continuous on $I$ (this implies that
$\{ \overline{A_\hbar} \}_{\hbar\in I}$ forms a continuous fields of
$C^*$-algebras over $I$, where $\overline{A_\hbar}$ is the
$C^*$-completion of $A_\hbar$), and \newline
(2) for $f,g\in A$,
$$\left\| {f*_\hbar g-g*_\hbar f \over {\sqrt{-1}\hbar}}-\{ f,g\}
\right\|_\hbar \to 0$$
as $\hbar\to 0$.  This definition still makes sense for a simple
Poisson
orbifold $M/\Gamma$, where $\Gamma$ is a finite group acting on
$M$ and the action preserves the Poisson bracket. A function on
$M/\Gamma$ is just a function on $M$ which is constant on each
orbit. A function $f$
on $M/\Gamma$ is defined to be smooth if it is smooth as a function on
$M$.
\enddefinition

\definition{Definition 2} A strict deformation quantization is {\it
flabby} if $A$ as above, contains $C^{\infty}_c(M)$, the algebra
of smooth functions of compact support on $M$. This notion also
makes sense for $M/\Gamma$ as above.
\enddefinition

There is more algebraic version of deformation quantization,
called {\it the formal deformation quantization}.
A formal deformation quantization of $M$ is defined as an
associative algebra structure $*$ on
$C^{\infty}(M)[[\hbar]]$ ( $\hbar$ is a formal letter ) such that,
 for $f, g \in C^{\infty} (M)$,
 $$f*g=fg+{\sqrt{-1}\over 2}\{ f,g\}\hbar+B_2(f,g)\hbar^2+B_3
 (f,g)\hbar^3+\cdots$$
 where $B_i$'s are bidifferential operators.
 Using ideas from string theory, Maxim Kontsevich proved that
 any Poisson manifold is formally deformation quantizable \Kont.

\smallskip

Now, we consider the example of a strict deformation quantization of
tori \RiefNon.
We use real coordinates $(x_1,\cdots,x_n)$ for the $n$-torus $T^n$,
viewing $T^n$ as $\Bbb{R}^n/\Bbb{Z}^n$. Any real
skew-symmetric matrix $\Theta$ defines a Poisson bracket on
$C^\infty(T^n)$:
$$\{ f, g\} := \sum_{j,k} \theta_{jk} {\partial f \over \partial x_j}
 {\partial g \over \partial x_k},$$
 for $f,g\in C^\infty(T^n)$.

The Fourier transform $\Cal{F}$  maps
$C^\infty(T^n)$ to $S(\Bbb{Z}^n)$, the space of complex-valued
Schwartz functions. Recall that $\Cal{F}(f)\in S(\Bbb{Z}^n)$, $f\in
C^\infty(T^n)$, is defined as follows:
for $p\in \Bbb{Z}^n$,
$$\widehat{f}(p):=\Cal{F}(f)(p)=
\int_{T^n} \exp(2\pi\sqrt{-1} x\cdot p)f(x)dx$$
where $dx$ is the Haar measure with $\int_{T^n}1dx=1$.
$\Cal F$ is invertible, and its inverse is given by
$$\phi\to \sum_{p\in \Bbb{Z}^n}\phi(p)\exp(2\pi \sqrt{-1} p\cdot x).$$
 $\Cal F$ carries the Poisson bracket to
$$\{\phi,\psi\}(p) =-4\pi^2\sum_{q}\phi(q)\psi(p-q)\gamma(q,p-q)$$
for $\phi,\psi\in S(\Bbb{Z}^n)$, where
$$\gamma(p,q)=\sum_{j,k}\theta_{jk}p_j q_k.$$

 For any $\hbar\in \Bbb{R}$, we define a function $\sigma_\hbar$
on $\Bbb{Z}^n\times\Bbb{Z}^n$ by
$$\sigma_\hbar(p,q)=\exp(-\pi \sqrt{-1} \hbar\gamma(p,q)),$$
and then define a deformed convolution product
$$(\phi *_\hbar\psi)(p)=\sum_q\phi(q)\psi(p-q)\sigma_\hbar(q,p-q).$$

The involution on $S(\Bbb{Z}^n)$ is defined, independent of $\hbar$,
as follow:
$$\phi^*(p)=\overline{\phi}(-p)$$
which, under the inverse of the Fourier transform, is just the
complex conjugation on $C^\infty (T^n)$.

Define a norm $\|\cdot\|_\hbar$ on $S(\Bbb{Z}^n)$ as the
operator norm for the action of $S(\Bbb{Z}^n)$ on $l^2(\Bbb{Z}^n)$
given by $\phi\cdot\xi=\phi *_\hbar\xi.$ We define $C_\hbar$ to be
$C^\infty(T^n)$ with the product, the involution, and the norm obtained by
pulling back, via the Fourier transform, the product
$*_\hbar$, the involution, and norm $\|\cdot\|_\hbar$ we defined above. Then
$\{ C_\hbar\}_{\hbar\in \Bbb{R}}$ is a strict deformation
quantization of the Poisson manifold $(T^n, \Theta)$.

$A_\Theta$ is defined as  $C_1$,  the algebra
for $\hbar=1$. Then, by definition, $C_\hbar =A_{\hbar\Theta}$.
 An easy computation shows that
$$U_k U_j=\exp(2\pi i\theta_{jk}) U_j U_k,$$
where $U_i=\exp(2\pi\sqrt{-1} x_i)$. The enveloping
$C^*$-algebra $\overline{A}_\Theta$ of $A_\Theta$
is the universal $C^*$-algebra generated by $n$
unitary operators satisfying the above relations.
$\overline{A}_\Theta$ is called {\it the noncommutative torus}.
The noncommutative tori appear naturally in M-theory
compactification. See \CDS.

We define a $\Bbb{Z}_2$-action on $T^n$ by
$$\gamma\cdot (x_1,\cdots,x_n)=(-x_1,\cdots, -x_n)$$
where $\gamma$ is the non-identity element of $\Bbb{Z}_2$.
(From now on, $\gamma$ will denote the non-identity element of
$\Bbb{Z}_2$.)
This $\Bbb{Z}_2$-action on $(T^n,\Theta)$ preserves the
Poisson bracket, that is,
$$\{ f^\gamma, g^\gamma \} =\{ f,g\}^\gamma$$
where $f^\gamma$ is defined as $f^\gamma(x)=f(-x)$. Also, the
strict deformation quantization of $(T^n,\Theta)$ defined as above is invariant
under the $\Bbb{Z}_2$-action, that is,
$$f^\gamma *_\hbar g^\gamma=(f*_\hbar g)^\gamma.$$
Hence the strict deformation quantization of $(T^n,\Theta)$
restricts to a strict deformation quantization of the Poisson
orbifold $T^n/\Bbb{Z}_2$, which is flabby. A smooth function $f$
on $T^n/\Bbb{Z}_2$ is just an smooth even function on $T^n$, that
is, $f(-x)=f(x).$
 This strict deformation
quantization is given by
$\{A_{\hbar\Theta}^\sigma\}_{\hbar\in\Bbb{R}}$.
$A_{\Theta}^\sigma$ denotes the subalgebra of
$A_{\Theta}$ which consists of even functions in $A_{\Theta}$. Its closure
$\overline{A}_{\Theta}^\sigma$ in
$\overline{A}_{\Theta}$ consists of even functions in
$\overline{A}_{\Theta}$. $\overline{A}_{\Theta}^\sigma$ is called
the symmetrized noncommutative torus.

\smallskip

We will simply write $U_p$ for $\exp(2\pi\sqrt{-1} p\cdot x)$,
$p\in \Bbb{Z}^n$.
Note that $(U_p)^\gamma=U_{-p}=(U_p)^*$. The  difference between
the action by $\gamma$ and the $*$-operation is that the former is
linear but the latter is conjugate-linear.
The algebra generated by $\{ U_p+U_{-p} \text{\ }
| \text{\ }p\in \Bbb{Z}^n \}$ is dense in
$\overline{A}_{\Theta}^\sigma$.

\proclaim{Theorem 1}
Assume that there exists an entry $\theta_{jk}$ of $\Theta$ such
that $4\theta_{jk}$ is not an integer. Then
the symmetrized noncommutative torus $\overline{A}_{\Theta}^\sigma$
is Morita-equivalent to $\overline{A}_\Theta\rtimes\Bbb{Z}_2$.
\endproclaim
\demo{Proof} (For the notion of Morita-equivalence, see
\RaWi.)
 We let $C$ and $D$ denote the algebra $C(\Bbb{Z}_2,A_\Theta)$ and the dense
 subalgebra $A_{\Theta}^\sigma$ of $\overline{A}_\Theta^\sigma$,
respectively, where $C(\Bbb{Z}_2,A_\Theta)$ is the set of maps
from $\Bbb{Z}_2$ to $A_\Theta$.
 Recall that the product on $C(\Bbb{Z}_2,A_\Theta)$ is given as follows:
 $$(\Lambda\Psi)(e)=\Lambda(e)\Psi(e)+\Lambda(\gamma)\Psi(\gamma))^\gamma$$
 and
 $$(\Lambda\Psi)(\gamma)=\Lambda(e)\Psi(\gamma)+\Lambda(\gamma)(\Psi(e))^\gamma,$$
where $e$ is the additive identity of $\Bbb Z_2$.
 For a $C$-$D$ bimodule, we take $\Cal{E}=A_\Theta$.
The right $D$-module structure on $\Cal{E}$ is given by right
multiplications. A $D$-valued
inner product on $\Cal{E}$ is defined by
$$\left < U,V \right >_D= U^* V+(U^*)^\gamma V^\gamma.$$
The left $C$-module structure on $\Cal{E}$ is given as follows:
for $\Psi\in C$, $U\in \Cal{E}$,
$$\Psi\cdot U =  \Psi(e)U+\Psi(\gamma)U^\gamma.$$
We define a $C$-valued inner product:
$$\left < U,V \right >_C(e)=UV^* \text{, \ \ } \left <U,V\right
>_C(\gamma)=U(V^*)^\gamma.$$
Then, easily, we have
$$\left < U,V \right >_C \cdot W=U\cdot \left <V,W\right >_D,$$
which is one of the requirements in the definition of
Morita-equivalence. \par
We proceed to prove that the linear span
$\left<\Cal{E},\Cal{E}\right>_C$ of
 $\{ <x,y>_C | x,y \in \Cal{E}
\}$ is all of $C$. Since $\left<\Cal{E},\Cal{E}\right>_C$
is not just a vector space but an ideal of $C$, we only need to
show that the identity element $\Phi_0$ of $C$ lies in
$\left<\Cal{E},\Cal{E}\right>_C$, where the identity element
$\Phi_0$ is given by $\Phi_0(e)=1$ and $\Phi_0(\gamma)=0$.
By the assumption, we have an entry
$\theta_{jk}$
such that $4\theta_{jk}$ is not an integer. We define an element
$\Lambda\in C$ by
$\Lambda(e)=U_j^{-2}$, $\Lambda(\gamma)=-U_j^{-2}U_k^2 U_j^2$.
Then we have
$$\Lambda \left ( \left< U_j, U_j^{-1}\right>_C-\left< U_k,
U_k^{-1}\right>_C+\left< U_k^2,1\right>_C
 \right )
 =(1-e^{8\pi\sqrt{-1}\theta_{jk}})\Phi_0.$$
Since $1-\exp(8\pi\sqrt{-1}\theta_{jk})$ is different from 0,
the identity element $\Phi_0$ lies in $\left< \Cal{E},\Cal{E}
\right>_C.$
Therefore $\left<\Cal{E},\Cal{E}\right>_C$ is dense in
$\overline{A}_\Theta\rtimes \Bbb{Z}_2$. \par
It is clear that
$\left<\Cal{E},\Cal{E}\right>_D$ is dense in
$\overline{A}_\Theta^\sigma$.
 Indeed, for any
$p\in \Bbb{Z}^n$,
$$\left < 1, U_{-p} \right>_D=U_p+U_{-p}.$$
 The inequalities
required in the definition of Morita-equivalence are also clearly
satisfied. Therefore, $\Cal E$ completes into a
Morita-equivalence bimodule between
$\overline{A}_\Theta\rtimes\Bbb{Z}_2$ and
$\overline{A}_\Theta^\sigma$.
$\square$
\enddemo

Hence $\overline{A}_\Theta\rtimes\Bbb{Z}_2$ and
$\overline{A}_\Theta^\sigma$ have the same $K$-groups,
provided $\Theta$ satisfies the assumption in the above theorem.
Therefore, we have the following theorem, which was proved for the
case of $\Theta$ being totally irrational \FaW.

\proclaim{Theorem 2} Assume that $\Theta$ has an entry
$\theta_{jk}$ such that $4\theta_{jk}$ is not an integer. Then
$$K_0(\overline{A}_\Theta^\sigma)=\Bbb{Z}^{3\cdot 2^{n-1}} \text{
and \ } K_1(\overline{A}_\Theta^\sigma)=0.$$
\endproclaim
\demo{Proof}
 Since
$K_0(\overline{A}_\Theta)=K_1(\overline{A}_\Theta)=\Bbb{Z}^{2^{n-1}}$
\RiefPro \
for any real skew-symmetric matrix $\Theta$, the same reasoning as
in Theorem~7 of
\FaW \  also works for this case. Hence, Theorem~1 implies this
theorem.
 $\square$
\enddemo

 For an abelian group $G$, we define $\operatorname{rk}(G)$ as the
 rank of the free abelian group $G/\operatorname{tor}(G)$. Here
 $\operatorname{tor}(G)$ denotes the torsion subgroup of $G$.

\proclaim{Theorem 3}
$\operatorname{rk}\left( K^0(T^4/\Bbb{Z}_2)\right )$ is greater than 24.
\endproclaim
\demo{Proof}
The space $T^4/\Bbb{Z}_2$ has 16 singularities, which we enumerate
by $p_1,  \cdots, p_{16}$.
 If we blow them
up, we obtain a $K3$ surface $Z$. (For $K3$ surfaces and their
significance in string theory, see \Asp.) We let $X_k$ be the
space obtained from $T^4/{\Bbb Z}_2$
by blowing up the first $k$ points $p_1, \cdots,
p_k$.  Hence $X_{16}=Z$. We consider the pair $(S^3, Z)$,
where $S^3$ is the 3-sphere to which the point $p_{16}$ has been blown up.
Then the one-point compactification of $Z-S^3$ is $X_{15}$. We
consider the following 6-term exact sequence in $K$-theory:
$$ \CD
K^0(Z-S^3) @> q^*_0 >>  K^0(Z)  @>i^*_0>>  K^0(S^3) \\
@A\partial_1AA  @.      @VV\partial_0V \\
K^1(S^3) @<i^*_1<< K^1(Z) @<q^*_1<< K^1(Z-S^3).
\endCD
$$
From the fact that $H^{ev}(Z,\Bbb{Q})\cong \Bbb{Q}^{24} $ and
$H^{odd}(Z,\Bbb{Q})=0$, we have
$$\operatorname{rk}(K^0(Z))=24 \text{, \ and }
\operatorname{rk}(K^1(Z))=0.$$
Since $\operatorname{im}i^*_1$ is a subgroup of $K^1(S^3)\cong
\Bbb{Z}$, $\operatorname{im}i^*_1$ must be $0$. It means that
$\partial_1$ is injective. Therefore $\operatorname{ker}q^*_0
\cong \Bbb{Z}$. Since $K^0(Z)/\operatorname{im}q^*_0$ is
isomorphic to a subgroup of $K^0(S^3)\cong \Bbb{Z}$, we have
$\operatorname{rk}(\operatorname{im}q^*_0)\geq \operatorname{rk}
(K^0(Z))-1=23$. The fact that
 $\operatorname{rk}(\operatorname{ker}q^*_0)=1$ implies that
$\operatorname{rk}(K^0(Z-S^3))\geq 24.$  Consequently, we have
$$\operatorname{rk}(K^0(X_{15}))=
\operatorname{rk}(K^0(Z-S^3)\oplus \Bbb{Z}) \geq 25.$$
Now, for $k=1,\cdots, n-1$, we consider the pair $(S^3, X_k)$,
where $S^3$ is the 3-sphere to which the point $p_k$ has been
blown up. Then the exat sequence
$$\CD
K^0(X_k-S^3) @>>> & K^0(X_k) @>>> & K^0(S^3) \cong \Bbb{Z}
\endCD
$$
gives us the inequality $\operatorname{rk}(K^0(X_k-S^3)) \geq
\operatorname{rk}(K^0(X_k))-1$. Since $X_{k-1}$ is the one-point
compactification of $X_k-S^3$, we have
$$\operatorname{rk}(K^0(X_{k-1})) \geq
\operatorname{rk}(K^0(X_k)).$$
Hence it follows that $\operatorname{rk}(K^0(T^4/\Bbb{Z}_2)) \geq
25.$  $\square$
\enddemo

\remark{Remark} Let $\Theta$ be
a nonzero $4\times 4$ real skew-symmetric matrix. Then, for a
generic number $s$, $s\Theta$ satisfies the assumption of
Theorem~1. Hence
$K_0(\overline{A}_{s\Theta}^\sigma)=\Bbb{Z}^{24}$, which is
different from
$K_0(\overline{A}_{0\cdot\Theta}^\sigma)=K_0\left(C(T^4/\Bbb{Z}_2)\right)
= K^0(T^4/\Bbb{Z}_2)$.
Therefore, the flabby strict deformation quantization $\{
A_{t\Theta}^\sigma \}_{t\in \Bbb{R}}$ of the Poisson orbifold
$T^4/\Bbb{Z}_2$
 gives us an example, where
$K_*(\overline{A}_{t\Theta}^\sigma)$ varies as $t$ varies.
\endremark

\Refs
\endRefs
\biblio

\enddocument

\end